\documentclass{scrartcl}
\usepackage[utf8]{inputenc}
\usepackage{amsfonts, amssymb, amsmath, amsthm, mathtools}
\usepackage[all]{xy}
\usepackage{tikz-cd}
\usepackage{soul} 
\usepackage[colorinlistoftodos]{todonotes}
\usepackage{tcolorbox}
\usepackage{quiver}

\usepackage{aliascnt}
\definecolor{greencolor}{rgb}{0,0.45,0}
\definecolor{ucphcolor}{rgb}{0.517,0.016,0.016}
\usepackage[colorlinks,
            citecolor=greencolor,
            linkcolor=ucphcolor,
            urlcolor=greencolor, unicode,linktocpage]{hyperref}
\usepackage[capitalize]{cleveref}

\usepackage[backend=biber, style=alphabetic, maxnames=10, maxalphanames=3, minalphanames=3]{biblatex}
\usepackage[margin = 1.25in]{geometry}

\usepackage{mathrsfs}
\usepackage{eucal}

\theoremstyle{definition}
\newtheorem{rmk}{Remark}[section]
\newtheorem*{war*}{Warning}
\newtheorem{dfn}[rmk]{Definition}
\newtheorem{defn}[rmk]{Definition}
\newtheorem{ex}[rmk]{Example}
\newtheorem{exm}[rmk]{Example}
\newtheorem*{ques*}{Question}
\newtheorem{goal}[rmk]{Goal}

\newtheorem*{ans*}{Answer}

\newtheorem{rec}[rmk]{Recollection}

\newtheorem{cons}[rmk]{Construction}
\newtheorem{cnstr}[rmk]{Construction}
\newtheorem{var}[rmk]{Variant}

\theoremstyle{plain}
\newtheorem{cor}[rmk]{Corollary}

\newtheorem{lem}[rmk]{Lemma}
\newtheorem{thm}[rmk]{Theorem}
\newtheorem*{thm*}{Theorem}
\newtheorem{prop}[rmk]{Proposition}

\graphicspath{{pics/}}

\def\A{\mathbb{A}}

\def\E{\mathbb{E}}
\def\F{\mathbb{F}}

\def\N{\mathbb{N}}

\def\Ss{\mathbb{S}}

\def\W{\mathbb{W}}
\def\Z{\mathbb{Z}}

\def\CC{\mathcal{C}}
\def\DD{\mathcal{D}}

\newcommand{\bb}[1]{\mathbb{#1}}

\newcommand{\mc}[1]{\mathcal{#1}}

\newcommand{\cochain}[2]{C^{*}(#1;#2)}
\newcommand{\DLM}{\rm{DLM}}

\newcommand{\rm}[1]{\mathrm{#1}}

\newcommand{\Mod}{\rm{Mod}}

\newcommand{\Map}{\mathrm{Map}}

\def\Spc{\mathcal{S}}
\def\Sp{\mathrm{Sp}}
\newcommand{\pushout}{\arrow[ul, phantom, "\ulcorner", very near start]}
\newcommand{\pullback}{\arrow[dr, phantom, "\lrcorner", very near start]}

\definecolor{DefColor}{rgb}{0.6,0.15,0.25}
\newcommand{\mdef}[1]{\textcolor{DefColor}{#1}}
\definecolor{twoColor}{rgb}{0.75,0,0.75}
\newcommand{\ptwo}[1]{\textcolor{twoColor}{#1}}

\newcommand{\psq}[2]{\beta^{#1}\!Q^{#2}}

\def\Aff{\mathrm{Aff}}
\def\DAff{\mathrm{DAff}}
\def\SAff{\Ss\mathrm{Aff}}

\def\E{\mathbb{E}}

\DeclareMathOperator{\Shv}{\mathrm{Shv}}
\DeclareMathOperator{\Alg}{\mathrm{Alg}}
\DeclareMathOperator{\colim}{\mathrm{colim}}

\DeclareMathOperator{\fib}{\mathrm{fib}}
\def\op{\mathrm{op}}

\def\CAlg{\mathrm{CAlg}}
\def\DAlg{\mathrm{DAlg}}
\def\DLA{\mathrm{Alg}_{\mathrm{\scriptscriptstyle DL}}}
\def\DLAP{\mathrm{Alg}_{\mathrm{\scriptscriptstyle DL}}^{\sharp}}
\def\DLM{\mathrm{Mod}_{\mathfrak{B}}^{i}}
\def\DLMP{\mathrm{Mod}_{\mathfrak{B}}^{i,\sharp}}
\def\LMod{\mathrm{LMod}}

\def\SymDL{\mathrm{Sym}_{\mathrm{\scriptscriptstyle DL}}}
\def\Sym{\mathrm{Sym}}

\allowdisplaybreaks

\tikzcdset{scale cd/.style={every label/.append style={scale=#1},
    cells={nodes={scale=#1}}}}

\title{Fixed points via tilting}
\author{Robert Burklund \and Vignesh Subramanian}
\date{}

\begin{document}

\maketitle

\begin{abstract}
    In this article we initiate the study of perfect commutative $\F_p$-algebras in the non-connective setting. As an application we show that given a finite $C_p$-space $X$ the $\F_p$-cochains on $X^{C_p}$ can be recovered as a base-change of the tilt of the Tate construction on cochains on $X^e$.
\end{abstract}

\tableofcontents
\vspace{5mm}


In this article we initiate the study of \mdef{perfect}
commutative $\F_p$-algebras. This being introduced in \ref{dfn:perfect} to single out those commutative $\F_p$-algebras $R$ whose frobenius Dyer--Lashof operation 
\[ Q^0 : \pi_*R \to \pi_*R \]
is an isomorphism.

Our motivating example is
the perfect algebra $\cochain{X}{k}$
of $k$-valued cochains on a space $X$
(here and throughout $k$ will be a perfect field of char. $p$). 
Indeed, after developing the basics of the theory of perfect algebras in \S\ref{sec:perfalg}
we recover several deep theorems in $p$-adic homotopy theory via straightforward arguments in \S\ref{sec:cochains}.

\subsection*{Perfect algebras}

Our study of perfect $\F_p$-algebras begins with establishing the following package of non-connective refinements of results familiar from the discrete setting:

\begin{prop} 
The inclusion of perfect $\F_p$-algebras
into all commutative $\F_p$-algebras admits
left and right adjoints:
\[ \begin{tikzcd} 
\CAlg_{\F_p}^{\sharp} \ar[r, hook, "\iota"] &
\CAlg_{\F_p} 
\ar[l, bend right, "(-)^{\sharp}"']
\ar[l, bend left, "(-)^{\flat}"]
\end{tikzcd} \]
perfection and tilting respectively.
This adjoint triple enjoys the following properties:
\begin{itemize}
\item Perfect $\F_p$-algebras are automatically coconnective.
\item The perfection $\F_p\{x_n\}^\sharp$ of a free $\F_p$-algebra on a class in degree $n < 0$ is well behaved in the sense that $\pi_*\F_p\{x_n\}^\sharp$ is a free graded commutative algebra on certain explicit generators in negative degrees.
\item The underlying spectrum of the tilt of a 
    commutative $\F_p$-algebra $R$ is the limit
    \[ R^\flat \to \cdots \xrightarrow{\mathfrak{Q}^0} R \xrightarrow{\mathfrak{Q}^0} R \xrightarrow{\mathfrak{Q}^0} R \]
    along a stable power operation that acts by $Q^0$ on homotopy groups.
\end{itemize}
\end{prop}

Note however that, unlike in the discrete setting, the frobenius power operation cannot be lifted to an endomorphism of the identity functor on $\CAlg_{\F_p}$.

A priori coconnectivity for perfect $\F_p$-algebras dovetails nicely with the theory of descent for bounded above modules over ring spectra (recalled in \S\ref{sec_coconndescent}). Using descent 
we then prove a unicity result for perfect geometry that we summarize with the slogan:

\begin{quote}
\textit{In the perfect setting spectral and derived algebraic geometry coincide.}
\end{quote}

\begin{thm}
    Let $\mathrm{DAlg}_{\F_p}$ denote the
    category of derived $\F_p$-algebras\footnote{%
    Meaning: See \cite{brantner2023deformation, raksit2020hochschild}.}
    and let $\mathrm{DAlg}_{\F_p}^\sharp$
    denote the full subcategory on the perfect derived $\F_p$-algebras.
    We construct isomorphisms
    \[ \CAlg_{\F_p}^{\sharp} \cong \mathrm{DAlg}_{\F_p}^{\sharp}, \]
    \[ \Shv_{\mathrm{hf}}^{\wedge}(\CAlg_{\F_p}^{\sharp,\op}) 
    \cong \Shv_{\mathrm{fpqc}}^{\wedge}(\CAlg_{\F_p}^{\heartsuit,\sharp,\op}) 
    \cong \Shv_{\mathrm{hf}}^{\wedge}(\mathrm{DAlg}_{\F_p}^{\sharp,\op})  \]
    of categories of perfect algebras
    and categories of accessible hypersheaves.
    Here the hf pretopology 
    is the non-connective 
    extension of the fpqc topology set up in \S\ref{sec_coconndescent}.
\end{thm}






Through this we classify nullstellensatzian\footnote{Meaning: See \cite{chroNS}.} perfect $\F_p$-algebras:

\begin{cor}
    The nullstellensatzians in $\CAlg_{\F_p}^{\sharp}$ are discrete algebraically closed fields.
\end{cor}

\subsection*{Algebraic Models of spaces}


Although the $k$-cochains functor may at first appear an unlikely candidate for an algebraic model of $p$-complete spaces, in a remarkable demonstration of the depths of complexity of commutative algebra structures, Mandell proved the following theorem:





\begin{thm}[{\cite{mandell2001algebras}}]
    For $k$ algebraically closed, the $k$-cochains functor
    \[ \cochain{-}{k}: \Spc \to \rm{CAlg}^{op}_k \]
    is fully faithful whenever the target is 
    a $p$-complete nilpotent space of finite type.
\end{thm}

In \S\ref{sec:MandellThm}, we 
offer a streamlined proof of Mandell's theorem showcasing the simplifying power of perfect $k$-algebras. Building on this we then prove a generic form of Yuan's extension \cite[Cor. 7.6.1]{yuan2023integral} of Mandell's theorem demonstrating that fully faithfulness of cochains is rather insensitive to the coefficient ring.

\begin{cor}
    Let $R$ be a $p$-complete connective commutative algebra 
    with $\pi_0R / p \cong k$ and alg. closed.
    The $R$-cochains functor is fully faithful
    whenever the target is the $p$-completion of a finite nilpotent space.
\end{cor}

\subsection*{Fixed points via Tilting}

Given a finite $C_p$-space $X$
its Euler characteristic satisfies
$\chi(X) \equiv \chi(X^{C_p}) \ \pmod{p}. $
Inspired in part by this elementary congruence, 
a major current in equivariant homotopy theory has been 
the study of \emph{Smith theory} (cf. \cite{smith1934theorem})
whose central goal we summarize as follows:

\begin{goal} \label{qst:motiviation}
    Let $X$ be a finite $C_p$-space.
    \begin{itemize}
        \item Recover the cohomology of $X^{C_p}$
        from the cohomology of $X^e$ with its $C_p$-action.
        \item Recover the space $X^{C_p}$ from $X^e$ with its $C_p$-action.
    \end{itemize}
\end{goal}



Landmark results in Smith theory include the localization theorems of Borel and Atiyah--Bott \cite{borel1960, atiyah1984moment} and the resolution of Sullivan's fixed point conjecture \cite{miller1984sullivan, carlsson1991equivariant, lannes1992espaces}. In \S\ref{sec:fp-tilt}, exploiting the observation that tilting allows us to de-localize localized cochains, we put cohomological Smith theory in near-optimal form:

\begin{thm}[Fixed points via Tilting]
    Let $X$ be a finite $C_p$-space.
    We construct a natural isomorphism of perfect $k$-algebras
    \[ \cochain{X^{C_p}}{k} \cong k \otimes_{k^{BC_p}} (\cochain{X^e}{k}^{tC_p})^{\flat}. \] 
\end{thm}

In view of Mandell's theorem, the theorem above is close in spirit to the Sullivan conjecture. Being inspired by this, in the final subsection we give a short proof that the Segal conjecture for $C_p$ \cite{lin1980conjectures, lin1980calculation, gunawardena1980segal} implies the nilpotent case of the Sullivan conjecture.

\section*{Acknowledgements}
The authors  would like to thank 
Tobias Bartel,
Natàlia Castellana,
Jesper Grodal, 
Gijs Heuts, 
Ishan Levy, 
Jesper Moller
for helpful conversations and comments on a draft.

During the course of this work
the authors were supported
by the Danish National Research Foundation
through the Copenhagen Center for Geometry and Topology (DNRF151).
The first author was additionally supported by Villum foundation through 
the Villum Young Investigator grant: Computational Frameworks in Homotopy Theory.
The second author was additionally supported by European Research Council (ERC) under Horizon Europe (GeoCats, grant No. 101042990), NWO VIDI Grant (no. 223.093) and would like to thank the Max Planck Institute for Mathematics (MPIM) in Bonn for their hospitality and excellent working conditions.

\section*{Notation and conventions}

\begin{enumerate}
\item We fix once and for all a prime $p$ and a perfect field $k$ of char. $p$.
\item If $p=2$ delete all \ptwo{purple} text.
\item If $\mc{C}$ is a stable category with a $t$-structure, we denote by $\mc{C}_{+}$, respectively $\mc{C}_{-}$, the full subcategories of bounded below, respectively bounded above, objects.
\item Let $\mc{S}, \mc{S}_*$ and $\rm{Sp}$ denote the category of spaces, pointed spaces and spectra, respectively.
    \item By a perfect $k$-algebra we mean a perfect commutative $k$-algebra (\ref{dfn:perfect}). 
    \item $k \{ x_n \}$ denotes the free commutative $k$-algebra on a generator $x_n$ in degree $n$.
    \item For brevity we let $\CAlg \coloneqq \CAlg(\Sp)$ and given an $R \in \CAlg$ we let
    $\CAlg_R \coloneqq \CAlg(\Mod_R) \cong \CAlg_{R/}$ and $\Map_R \coloneqq \Map_{\CAlg_R}$.
\end{enumerate}

\section{Perfect algebras} \label{sec:perfalg}

In this section, we develop the theory of perfect commutative $\F_p$-algebras. For us, an $A \in \CAlg_{\F_p}$ is perfect as soon as the frobenius Dyer--Lashof operation \[ Q^0 : \pi_*A \to \pi_*A \] is an isomorphism (see \ref{dfn:perfect}). We use $\sharp$ to indicate a restriction to perfect algebras. The technical heart of the section is \S\ref{sec:perf-tilt}. Therein we establish the following package of non-connective refinements of results familiar from the discrete setting:

\begin{itemize}
\item Perfect $\F_p$-algebras are automatically coconnective.
\item $\CAlg^{\sharp}_{\F_p} \hookrightarrow \CAlg_{\F_p}$ admits left and right adjoints: perfection and tilting. 
\item Perfections of free algebras are well-behaved.
\item The tilt can be computed as an inverse limit along $Q^0$.
\end{itemize}

A priori coconnectivity for perfect $\F_p$-algebras dovetails nicely with the theory of descent for bounded above modules over ring spectra (recalled in \S\ref{sec_coconndescent}). Using descent we embed perfect $\F_p$-algebras into fpqc (hyper)sheaves on \emph{discrete} perfect $\F_p$-algebras (\S\ref{sec:stacks}) and observe that spectral and derived perfect $\F_p$-algebras coincide (\S\ref{sec:derived}). Finally, we show that bounded below perfect $\F_p$-algebras admit unique spherical lifts (\S\ref{sec:sphericalwitt}).

\subsection{Dyer--Lashof operations}
\label{sec:DL}

Before beginning our analysis of perfect algebras it will be convenient to review the theory of homotopy power operations for commutative $\F_p$-algebras, i.e. Dyer--Lashof operations. The original reference for much of this material is \cite{dyer1962homology} and a more detailed exposition can be found across \cite{bruner1986h}, \cite{lawson2020en} and \cite{glasman2020stable}.

Let $\mathfrak{B}$ denote the graded $\F_p$-algebra of DL operations:
\begin{itemize}
\item\textbf{Generators}: $\psq{\epsilon}{i}$ of degree $\ptwo{2(p-1)}i - \epsilon$
    where $\epsilon \in \{0\ptwo{, 1}\}$, $i \in \Z$.\footnote{We take the convention that $\psq{\epsilon}{i} = 0$ if $\epsilon \not\in \{0\ptwo{,1}\}$ and drop the $\beta$ if its index is $0$.}
\item\textbf{Adem relations}: Let $c_{ij}^{k\ell} \coloneqq (-1)^{i+k} \binom{(p-1)(k-j) - 1 + \ell}{pk-i}$, then
    \begin{align*}
        \psq{\epsilon}{i} \cdot \psq{e}{j} 
        &= \sum_{\substack{\frac{i+s}{p} \leq k < i+e-(p-1)j \\ 0 \leq s \leq e}} c_{i+s,j}^{k,e-s} \psq{\epsilon+e-s}{i+j-k} \cdot \psq{s}{k}, & &i+e > pj.
    \end{align*}
\item\textbf{Basis}: Given 
    $I = (\epsilon_1,i_1,\dots,\epsilon_n,i_n) \in \left(\{0\ptwo{,1}\} \times \Z \right)^{n}$
    we let $Q^{I} \coloneqq \psq{\epsilon_1}{i_1} \cdots \psq{\epsilon_n}{i_n}$.
    \begin{itemize}
    \item We say $I$ is \mdef{admissible} if
    $\Delta_k \coloneqq p i_{k} - \epsilon_{k} - i_{k-1} \geq 0$ for $1 < k \leq n$.
    \item The admissible $Q^I$ are a basis for $\mathfrak{B}$.
    \item The \mdef{excess} of $I$ is 
    $e(I) \coloneqq 2i_n - \sum_{k=1}^n \epsilon_k - 2 \sum_{k=2}^n \Delta_k $.
    \item $I$ is \mdef{positive} (\mdef{negative}) if each term in it has positive (negative) degree.
    \item For $I$ admissible $I = J 0^k L$ with $J$, $L$ admissible, $J$ negative and $L$ positive.
    \end{itemize}
\end{itemize}

A \mdef{DL module} is a graded left $\mathfrak{B}$-module satisfying:
\begin{itemize}
    \item\textbf{Instability}: $\psq{\epsilon}{i} \cdot x = 0$ for $|\psq{\epsilon}{i}| < (p-1)|x|$.
\end{itemize}
We let $\DLM$ denote the category of DL modules 
and equip it with the symmetric monoidal structure 
created by the
\begin{itemize}
    \item\textbf{Cartan formula}:
    $ \psq{\epsilon}{i} \cdot (x \otimes y) = \sum_{\substack{j+k = i \\ a+b=\epsilon}} \left( \psq{a}{j} \cdot x \right) \otimes \left( \psq{b}{k} \cdot y \right) $.
\end{itemize}

Let $\mathfrak{B}\langle n \rangle$ denote 
the free DL module on a class $x_n$ in degree $n$.
This DL module has a basis of classes 
$Q^I \cdot x_n$ with $I$ admissible, $e(I) \geq n$.   
Let $q^I : \mathfrak{B}\langle m \rangle \to \mathfrak{B}\langle n \rangle$
denote the map sending $x_m$ to $Q^I \cdot x_n$ (here $m=n+|I|$).
The map $q^I$ co-represents the action of $Q^I$ from degree $m$ to degree $n$ on DL modules.

A \mdef{DL algebra} is a commutative algebra in $\DLM$ satisfying:
\begin{itemize}
    \item\textbf{Power compatibility}: $Q^k(x) = x^{p}$ when $\ptwo{2}k = | x |$.
\end{itemize} 
We denote the category of DL algebras by $\DLA$.
The underlying DL module functor admits a left adjoint $\SymDL$.
This sends a DL module $M$ to the maximal quotient of $\Sym(M)$
satisfying the power relation (the other relations being automatic).

\begin{rec}\label{rec: DLA}
The functor $\pi_*$ lifts to a functor
\[ \pi_\star : \CAlg_{\F_p} \to \DLA \]
enjoying the following properties:
\begin{enumerate}
\item $\pi_\star$ preserves coproducts and filtered colimits.
\item $\pi_\star \F_p\{x_n\} \cong \SymDL(\mathfrak{B}\langle n \rangle)$.
\item Let $\mathcal{Q}^I : \F_p\{ x_{m} \} \to \F_p\{ x_n \} $
    denote the map sending $x_m$ to $Q^I(x_n)$ (here $m=|I|+n$).
    Under (2), $\pi_\star \mathcal{Q}^I \cong \SymDL(q^I) $,
    i.e. $\pi_\star \mathcal{Q}^I$ co-represents the action of $Q^I$ from degree $n$ to degree $m$ on DL algebras.
\item The DL operations are stable in the weak sense that
    given a pullback square in $\CAlg_{\F_p}$ the associated long exact sequence
    on $\pi_*$ lifts to $\DLM$.
\item The underlying functor $\CAlg_{\F_p} \to \rm{Sp}$ is 
    co-represented by the co-spectrum object 
    \[ u := \left\{   \F_p\{ x_n\}, \qquad \begin{tikzcd}[sep=tiny] \F_p\{ x_n\} \ar[r] \ar[d] & \F_p \ar[d] \\ \F_p \ar[r] & \F_p\{ x_{n+1} \} \pushout \end{tikzcd}  \right\}_{n \in \bb{Z}}  \]
    in $\CAlg_{\F_p}$.
    DL operations are stable in the strong sense that 
    the maps $\mathcal{Q}^I$ assemble into  
    a map $ \mathfrak{Q}^I : \Sigma^{|I|} u \to u $
    acting by $Q^I$ on homotopy groups.
\end{enumerate} 
\end{rec}



\subsection{Perfection and tilting}
\label{sec:perf-tilt}

Building on our understanding of Dyer--Lashof operations, 
we define and study perfect $\F_p$-algebras.
The titular, perfection and tilting, are the left and right adjoints to
the inclusion of perfect $\F_p$-algebras into all commutative $\F_p$-algebras.

\begin{dfn} \label{dfn:perfect}
    Let $\mathfrak{B}^\sharp$ denote the localization of the Dyer--Lashof algebra at $Q^0$.\footnote{%
    Meaning: The universal graded ring under $\mathfrak{B}$ in which $Q^0$ is a two-sided unit.}
    Let $ \DLMP$ denote the full subcategory of $Q^0$-local DL modules and let the pullbacks
    \[ \begin{tikzcd}
        {\CAlg_{\F_p}^{\sharp}} \ar[r, "\pi_\star"] \ar[d, hook] \pullback &
        \DLAP \ar[r] \ar[d, hook] \pullback &
        \DLMP \ar[d, hook] \\
        {\CAlg_{\F_p}} \ar[r, "\pi_\star"] &
        \DLA \ar[r] &
        \DLM
    \end{tikzcd} \]
    create the corresponding full subcategories of 
    \mdef{perfect DL algebras} and \mdef{perfect $\F_p$-algebras}.
    We will use the terms perfect and $Q^0$-local interchangeably.
\end{dfn}

\begin{ex} \label{ex:disc-perfect}
    A discrete $\F_p$-algebra is perfect if its frobenius is an isomorphism.
\end{ex}

In our study of perfect $\F_p$-algebras we will need to understand $\mathfrak{B}^{\sharp}$.
The requisite control over this, a priori wild, localization of associative algebras comes from the Ore condition.

\begin{rec}[{\cite[Section~7.2.3]{lurieHA}}] \label{criteriaOrecondition}
    Let $R$ be a graded associative ring and 
    let $S \subset R$ be a set of homogeneous elements in $R$. 
    We say $S$ satisfies the \mdef{right Ore condition}  
    if the following conditions hold:
    \begin{enumerate}
        \item $S$ is multiplicatively closed.
        \item For every $r \in R$ and $ s \in S$ there exist $r' \in R$ and $s' \in S$ such that $rs' = sr'$.
        \item $x \in R$ and $s \in S$ such that $sx = 0$, there exist $s'$ such that $xs' = 0$.
    \end{enumerate} 
    If $S$ satisfies the right Ore condition, 
    then the underlying left $R$-module of
    the $S$-localization of $R$
    is given by the colimit along right multiplication by elements of $S$.
\end{rec}

\begin{var} \label{rmk:ore-var}
    Observe that given an $s \in R$ and a $T \subseteq R$ such that
    \begin{enumerate}
        \item $T$ generates $R$ and
        \item for each $t \in T$ we have $ts^n = sts^{n-1}$ for $n \gg 0$, 
    \end{enumerate}
    $\{ s^n \}_{n \geq 0}$ satisfies the right Ore condition and 
    $s$ is central $R[s^{-1}]$.
\end{var}


\begin{lem} \label{DyerLashofOre} \label{q0torsion}
    In $\mathfrak{B}$ we have 
    \begin{align}
        \psq{\epsilon}{i} \cdot (Q^0)^i &= 0 & &\text{ for }\quad i > 0,  \label{eq:oreoddPosdegree}\\
        \psq{\epsilon}{i} \cdot (Q^{0})^{1-(p-1)i} &= Q^{0} \cdot \psq{\epsilon}{i} \cdot (Q^{0})^{-(p-1)i} & &\text{ for }\quad i \leq 0. \label{eq:oreoddNegdegree}
    \end{align}
\end{lem}

\begin{proof}
    We prove \eqref{eq:oreoddPosdegree} by induction on $i > 0$ using the Adem relations:
    \[ \psq{\epsilon}{i} \cdot (Q^0)^i = \left( \sum_{\frac{i}{p} \leq k < i} c_{i,0}^{k,0} \psq{\epsilon}{i-k} \cdot Q^{k} \right) \cdot (Q^0)^{i-1} = 0. \]
    Let $i \leq 0$ and $r = -(p-1)i$. We prove \eqref{eq:oreoddNegdegree} using \eqref{eq:oreoddPosdegree} and the Adem relations:
    \[ Q^0 \cdot \psq{\epsilon}{i} \cdot (Q^0)^{r}  
        = \left( \sum_{\substack{\frac{s}{p} \leq k < \epsilon + r \\ 0 \leq s \leq \epsilon}} c_{s,i}^{k,\epsilon - s} \psq{\epsilon-s}{i-k} \cdot \psq{s}{k} \right) \cdot  (Q^0)^{r}  
        = \psq{\epsilon}{i} \cdot (Q^{0})^{1+r}.\hfill\qedhere \]
\end{proof}

    


\begin{prop} \label{cor:perfect-props}
$\mathfrak{B}^\sharp$ and the perfect DL modules enjoy the following properties:
\begin{enumerate}
    \item $Q^0$ is central in $\mathfrak{B}^\sharp$.
    \item $\mathfrak{B}^{\sharp}$ has a basis $\{ Q^J (Q^0)^{k} \}$ 
    with $J$ admissible, negative, $k \in \Z$. 
    \item The kernel of $\mathfrak{B} \to \mathfrak{B}^\sharp$ has a basis consisting of admissible classes $Q^I$ whose final index is positive.
    \item Every perfect DL module is concentrated in non-positive degrees.
    \item The perfect DL modules are closed under limits and colimits.
    \item The perfect DL modules are closed under tensor products.
    \item The free perfect DL module on a generator $x_n$ in degree $n$ is given by 
    \[ \mathfrak{B}\langle n \rangle^{\sharp} \coloneqq \colim \left( \mathfrak{B}\langle n \rangle \xrightarrow{q^0} \mathfrak{B}\langle n \rangle \xrightarrow{q^0} \cdots \right). \]
    This has a basis $\{ Q^J (Q^0)^{k} \cdot x_n \}$ 
    with $e(J) \geq n$, $J$ admissible, negative, $k \in \bb{Z}$.
\end{enumerate}
\end{prop}

\begin{proof}
    \ref{q0torsion} verifies the the hypotheses of variant \ref{rmk:ore-var} of \ref{criteriaOrecondition}.
    Using \ref{criteriaOrecondition} and \ref{q0torsion} we obtain (1), (2) and (3). 
    For (4) we observe that instability tells us the isomorphism $Q^0$ acts by zero in positive degrees.
    (5) is a general feature of localizations.
    For (6), comparing the Cartan formula with the vanishing from (3)
    we find that $Q^0$ acts on $M \otimes N$ by the isomorphism $Q^0 \otimes Q^0$.

    Comparing \ref{q0torsion}(\ref{eq:oreoddPosdegree}) with the basis of $\mathfrak{B}\langle n \rangle$ 
    we obtain the desired basis for $\mathfrak{B}\langle n \rangle^{\sharp}$.
    As $q^0$ co-represents $Q^0$ in degree $n$
    and this restricts to an isomorphism on perfect DL modules,
    proving the universal property of $\mathfrak{B}\langle n \rangle^{\sharp}$ reduces to knowing that
    $\mathfrak{B}\langle n \rangle^{\sharp}$ is perfect.
    Comparing the basis for $\mathfrak{B}\langle n \rangle^{\sharp}$ with the relations from \ref{q0torsion}(\ref{eq:oreoddNegdegree}) 
    we learn that $\mathfrak{B}\langle n \rangle^{\sharp}$ is indeed perfect.
\end{proof}


\begin{lem} \label{lem:perfect-props}
    $\CAlg_{\F_p}^{\sharp} \xhookrightarrow{\iota} \CAlg_{\F_p}$
    is closed under (co)limits.
    Invoking an adjoint functor theorem, 
    $\CAlg_{\F_p}^{\sharp}$ is then presentable and $\iota$
    admits left and right adjoints.
    \begin{itemize}
        \item We write \mdef{$(-)^{\sharp}$} for the left adjoint of $\iota$ and refer to it as \mdef{perfection}.
        \item We write \mdef{$(-)^{\flat}$} for the right adjoint of $\iota$ and refer to it as \mdef{tilting}.
    \end{itemize}
\end{lem}

\begin{proof}
    An $R \in \CAlg_{\F_p}$ is perfect if and only if $\mathfrak{Q}^0(R)$ is an isomorphism.
    As limits and colimits of isomorphisms are isomorphisms and 
    $u$ commutes with limits and sifted colimits, 
    we conclude that perfect $\F_p$-algebras are closed under limits and sifted colimits.
    For binary coproducts compare \ref{rec: DLA}(1) and \ref{cor:perfect-props}(6).
    The empty colimit gives $\F_p$ and
    this is perfect (\ref{ex:disc-perfect}).
\end{proof}

\begin{ex} \label{ex:cochainperf}
  Let $R$ be a discrete perfect $\F_p$-algebra.
  As perfect $\F_p$-algebras are closed under limits, 
  the cochain algebra $C^*(X;R)$ is perfect for all spaces $X$.
\end{ex}

\begin{prop} \label{thm:free-perfect}
    We construct isomorphisms
    \begin{align*}
       \F_p\{x_n\}^{\sharp} 
       &\cong \colim \left( \F_p\{x_n\} \xrightarrow{\mathcal{Q}^0}
       \F_p\{x_n\} \xrightarrow{\mathcal{Q}^0} \cdots \right), \\
       \pi_\star \F_p\{x_n\}^\sharp 
       &\cong \SymDL( \mathfrak{B}\langle n \rangle^{\sharp} ).
    \end{align*}
\end{prop}

\begin{proof}
    Let $A$ be the colimit on the rhs of the first claim.
    As $\mathcal{Q}^0 \cong \Omega^{\infty + n}\mathfrak{Q}^0$ 
    and this restricts to an isomorphism on perfect $\F_p$-algebras (\ref{rec: DLA}(5))
    the first claim reduces to perfectness of $A$.
    Using \ref{rec: DLA}(1,2,3) we calculate that
    \begin{align*}
        \pi_\star(A)
        \cong \pi_\star(\colim_{\mathcal{Q}^0} \F_p\{x_n\}) 
        \cong \SymDL( \colim_{q^0} \mathfrak{B}\langle n \rangle ) 
        \cong \SymDL(\mathfrak{B}\langle n \rangle^{\sharp} ).
    \end{align*}
    Observe that $\mathfrak{B}\langle n \rangle^{\sharp}$ is perfect (\ref{cor:perfect-props}(7)) and 
    $\SymDL$ preserves perfectness (\ref{cor:perfect-props}(5,6)).
\end{proof}

\begin{cor} \label{homotopygroupsoftilt}
    Given $A \in \CAlg_{\F_p}$,
    the underlying spectrum of $A^\flat$ is 
    the limit along $\mathfrak{Q}^0$ of the underlying spectrum of $A$.
    In particular, we have short exact sequences
    \[ 0 \to \lim_{Q^0}\!{}^{1}\ \pi_{n+1}(A)\to \pi_{n}( A^{\flat}) \to \lim_{Q^0} \pi_{n}(A) \to 0. \]
\end{cor}

\begin{proof}
    Let $u^{\sharp}$ denote the colimit of $u$ along $\mathfrak{Q}^0$.
    Using \ref{thm:free-perfect} we calculate that
    \[ u(A^{\flat}) \cong u^{\sharp}(A^{\flat}) \cong u^{\sharp}(A) \cong \lim_{\mathfrak{Q}^0} u(A).\hfill\qedhere \]
\end{proof}

\begin{rmk} \label{rmk:tilt-simplify}
    Perfect $\F_p$-algebras being coconnective, \ref{homotopygroupsoftilt}
    is only useful for $n \leq 0$.
    In the $n=0$ case, by instability and power compatibility
    $Q^0$ acts by zero on $\pi_1A$ and $x \mapsto x^p$ on $\pi_0A$
    so that $\pi_0A^{\flat}$ 
    is isomorphic to the classical tilt of $\pi_0A$.
\end{rmk}

\subsection{Descent} \label{sec_coconndescent}

In this interlude we analyze descent for bounded above modules over ring spectra. 
Our aim here is to isolate consequences of the standard proof of faithfully flat descent.
This material is comparatively elementary, but plays a key role in what follows.

\begin{defn} \label{defn:hfcover}
    Given a map of $\E_1$-algebras $f : A \to B$ 
    we let $I(f) \coloneqq \fib(A \to B)$
    and say that $f$ is a \mdef{half flat (hf) cover} if $I(f)$ has tor-amplitude $<0$ as a right $A$-module.\footnote{Meaning: For $M \in \LMod_A$ coconnective, $I(f) \otimes_A M$ is $(-1)$-coconnective.} 
\end{defn}

\begin{exm} \label{exm:hf-is-fpqc}
    If $A$ and $B$ are discrete commutative rings, then $A \to B$ is a hf-cover if and only if $A \to B$ is faithfully flat.
\end{exm}

\begin{lem} \label{exm:hf-cover-easy}
    Let $f : B \to C$ be a map of coconnective $\E_1$-algebras 
    with $A \coloneqq \pi_0B$ von Neumann regular%
    \footnote{Meaning: Every discrete left $A$-module is flat.}
    and $A \to \pi_0C$ injective. 
    The map $B \to C$ is an hf-cover.
\end{lem}

\begin{proof}
    Let $M \in \LMod(B)_{\leq 0}$.
    Expanding $I(f) \otimes_B M$ as the geometric realization of 
    $ I(f) \otimes_{A} B^{\otimes_A \bullet} \otimes_{A} M  $
    the bar filtration 
    (filter $\Delta^\op$ by $\Delta_{\leq *}^\op$)
    has associated graded
    \begin{equation} \label{eqn:bar-gr}
         I(f) \otimes_A (\Sigma \tau_{<0}B)^{\otimes_A * } \otimes_A M. 
    \end{equation}
    Our hypothesis that $A$ is von Neumann regular
    implies that $k$-coconnective $A$-modules
    have tor-amplitude $\leq k$.
    In particular, each term in (\ref{eqn:bar-gr}) is $(-1)$-coconnective.
    To conclude we observe coconnectivity is preserved by extensions and filtered colimits.
\end{proof}

\begin{lem} \label{lem:hf-props}
    The collection of hf-covers enjoys the following properties:
    \begin{enumerate}
    \item The hf-covers are closed under composition.
    \item The hf-covers are closed under filtered colimits in $\Alg^{\Delta^1}$.
    \item Given a square of $\E_1$-algebras 
    with $B \otimes_A C \to D$ an isomorphism,
    \[ \begin{tikzcd} 
    A \ar[r, "f"] \ar[d] & B \ar[d] \\
    C \ar[r, "g"] & D 
    \end{tikzcd} \] 
    if $f$ is an hf-cover, then $g$ is also an hf-cover.
    \end{enumerate}
\end{lem}

\begin{proof}
    \textit{Part (1)}
    Let $A \xrightarrow{f} B \xrightarrow{g} C \in \Alg$ be hf covers and 
    let $M \in \LMod_A$ be coconnective.
    Examining the cofiber sequences and isomorphism
    \begin{align*}
        &\left( I(f)  \to I(f \cdot g) \to I(g) \right) \otimes_A M  & 
        &I(g) \otimes_A M \cong I(g) \otimes_B (B \otimes_A M) \\
        &\left( A \to B \to \Sigma I(f) \right) \otimes_A M
    \end{align*}
    we conclude that $I(f \cdot g) \otimes_A M$ is $(-1)$-coconnective.


    
    
    \noindent\textit{Part (2).}
    Let $(f_\bullet : A_\bullet \to B_\bullet) : K \to \Alg^{[1]}$ 
    be a filtered diagram of $\E_1$-algebras consisting of hf-covers pointwise and
    let $f : A \to B$ denote its colimit.
    The isomorphism
    \[ I(f) \cong \colim_{k \in K} I(f_k) \otimes_{A_k} A \]
    witnesses $I(f)$ as a filtered colimit of right $A$-modules of
    tor-amplitude $\leq -1$.
    
    
    \noindent\textit{Part (3).}
    Using base-change we obtain isomorphisms
    \begin{align*}
        I(g) \otimes_C (-) 
        &\cong \fib(C \to B \otimes_A C) \otimes_C (-) \\
        &\cong \fib(A \to B) \otimes_A C \otimes_C (-) \\
        &\cong I(f) \otimes_A (-).\hfill\qedhere
    \end{align*} 
\end{proof}

\begin{cnstr} \label{cnstr:hf-pretopology}
    Through \ref{lem:hf-props}, 
    we may equip $\CAlg^{\op}$ with the \mdef{hf pretopology}:
    the finitary superextensive pretopology generated from the hf covers.
    Using this pretopology we instantiate several sites. 
    \begin{center}
    \begin{tabular}{|c|c|c|c|}\hline
        Site & Objects & Adjectives & Pretopology \\\hline\hline
        $\SAff$ & $\CAlg^{\op}$ & $\emptyset$ & hf \\\hline
        $\SAff_{-}$ & $\CAlg_{-}^{\op}$ & bounded above & hf \\\hline
        $\SAff_R$ & $\CAlg_R^{\op}$ & $R$-algebra & hf \\\hline
        $\SAff_{\F_p}^{\sharp}$ & $\CAlg_{\F_p}^{\sharp,\op}$ & perfect $\F_p$-algebra & hf \\\hline
        $\Aff_{R}$ & $\CAlg_{R}^{\heartsuit,\op}$ & discrete $R$-algebra & fpqc $=$ hf \\\hline
        $\Aff_{\F_p}^{\sharp}$ & $\CAlg_{\F_p}^{\sharp,\heartsuit,\op}$ & discrete, perfect $\F_p$-algebra & fpqc $=$ hf \\\hline
    \end{tabular}    
    \end{center}
\end{cnstr}

\vspace{6pt}

\begin{prop} \label{cor:modules-sheaf}
    Let $f : A \to B$ be a hf-cover of $\E_1$-algebras. The induced adjunction 
    \[ \widehat{f}^* : \widehat{\Mod}_A \rightleftharpoons \widehat{\Mod}_B : \widehat{f}_* \]
    on Postnikov completions of module categories is comonadic.
    As consequences,
    $\widehat{\Mod}$ is a hf-hypersheaf on $\SAff$ and
    the hypercompletion of hf is subcanonical on $\SAff_-$.
\end{prop}

\begin{proof}
    As upper-stars preserve connectivity we may write $\widehat{\Mod}$ as
    the limit of $\Mod^{\leq b} : \Alg \to \Pr^L$ along truncation maps
    (thereby verifying that the object under consideration is well-defined).
    Observe now that our claims live at the level of the wide subcategory of hf-covers.
    With this restriction upper-stars become $t$-exact so that 
    we may write $\widehat{\Mod}$ as the limit of $\Mod^{[a,b]}$.
    Finally, observe that 
    the inputs to Barr--Beck--Lurie are preserved under limits of categories, 
    completing our reduction to the truncated level.
    
    As $f^*$ is $t$-exact conservativity reduces to detecting zero objects.
    Given $M \in \Mod_A^{[a,b]}$ with $f^*M = 0$
    examining the induced isomorphism $\tau_{\geq a}(I(f) \otimes_A M) \cong M$
    we conclude by our tor-amplitude hypothesis on $I(f)$
    that $M$ has coconnectivity at most one less than itself, i.e. $M=0$.
    At the level of $[a,b]$-truncated objects $f^*$ preserves finite limits.
    On the other hand, at this level totalizations are finite limits.
    Thus, comonadicity.
    
    Consequence (1).
    As above we restrict to the wide subcategory on hf-covers.
    Through base-change isomorphisms comonadicity translates into 
    $\Mod^{[a,b]}$ being a hf-sheaf. 
    Truncated sheaves are automatically hypersheaves.
    Now pass to the limit over $a$ and $b$.
    
    Consequence (2).
    End of the unit in $\widehat{\Mod}$
    is a hf-hypersheaf on $\SAff$ whose restriction to 
    $\SAff_-$ is the underlying spectrum functor (conservative, preserves limits).   
\end{proof}

\subsection{Perfect stacks} \label{sec:stacks}

As an application of the theory developed up to this point,
we next show that perfect $\F_p$-algebras are in fact 
simply special examples of fpqc stacks on \emph{discrete} perfect $\F_p$-algebras.
Throughout this subsection we will use the results of \S\ref{sec:perf-tilt} freely.

\begin{thm} \label{thm:perfect-embed-discrete}
    The map of sites $\iota^\heartsuit : \Aff_{\F_p}^\sharp \to \SAff_{\F_p}^\sharp$ induces an isomorphism
    \[ \Shv(\Aff_{\F_p}^\sharp) \cong \Shv(\SAff_{\F_p}^\sharp) \]
    on categories of accessible hypersheaves 
    and through \ref{cor:modules-sheaf} an embedding of perfect $\F_p$-algebras
    into accessible fpqc hypersheaves on discrete perfect $\F_p$-algebras.
\end{thm}

\begin{lem} \label{lem:topos-basis}
    Let $i : \CC \hookrightarrow \DD$ be an embedding of op-accessible sites such that:
    \begin{enumerate}
        \item $\CC$ and $\DD$ are finitary and superextensive.
        \item The pretopology on $\CC$ is restricted from $\DD$.
        \item Every object of $\DD$ admits a cover by an object of $\CC$. 
    \end{enumerate}
    The induced map on categories of accessible hypersheaves is an isomorphism.
\end{lem}

\begin{proof}
    Consider the induced adjunction on topoi of hypersheaves
    \[ i^* : \Shv(\CC) \rightleftharpoons \Shv(\DD) : i_*. \]
    Using (1) we may re-express the hypersheaf condition on $\CC$ and $\DD$
    as product preserving plus sends semicosimplicial hypercovers to limits. 
    Using (2) and (3) every such hypercovering in $\DD$ is refined by one lying in $\CC$.
    From this we learn that 
    sheafification on $\DD$ commutes with restriction to $\CC$ and
    $i_*$ preserves colimits.    
    This lets us verify that $i_*$ is a colocalization by calculating that for $c \in \CC$ we have
    \begin{align*}
        i_*i^* y_\CC(c) 
        &\cong i_*y_{\DD}(i(c)) 
        \cong i_* L_\DD \Map_{\DD}(-,i(c)) 
        \cong L_{\CC} \Map_{\DD}(i(-),i(c)) 
        \cong y_{\CC}(c).
    \end{align*} 
    On the other hand,
    given $d \in \DD$ we may pick a hypercovering of $d$ lying entirely in $\CC$
    through which we witness that $y_{\DD}(d)$ lies in the image of $i^*$.
    As the $y_{\DD}(d)$ are generators of $\Shv(\DD)$ this implies $i_*$ is also conservative, 
    completing the proof.
\end{proof}

\begin{var} \label{var:small-object}
    Quillen's small object argument allows us to verify \ref{lem:topos-basis}(3) via:
    \begin{enumerate}
    \item[(3')] $\DD$ admits cofiltered limits and its singleton covers are closed under cofiltered limits.
    \item[(3'')] There exists a set of covers $\{ f_i : Y_i \to X_i \}_{i \in I}$ in $\DD$ such that
            any $Z \in \DD$ for which every map $Z \to X_i$ lifts through $Y_i$ for all $i\in I$ in fact lies in $\CC$.
    \end{enumerate}
\end{var}


\begin{proof}[Proof (of \ref{thm:perfect-embed-discrete}).]
    We prove the theorem by applying variant \ref{var:small-object} of \ref{lem:topos-basis}.
    The hf pretopology on perfect $\F_p$-algebras restricts to 
    the fpqc topology on discrete perfect $\F_p$-algebras (\ref{exm:hf-is-fpqc})
    with both pretopologies being finitary and superextensive by construction.
    Cofiltered limits of hf singleton hf covers are hf covers (\ref{lem:hf-props}(2)).
    The maps $\{  \F_p\{x_n\}^{\sharp} \to \F_p  \}_{n < 0}$
    each have coconnective source and induce an iso on $\pi_0$ (\ref{thm:free-perfect})
    so that they are hf covers by \ref{exm:hf-cover-easy}.
    Weak orthogonality to this set of maps is equivalent to connectivity,
    but recall that perfect $\F_p$-algebras are always coconnective (\ref{cor:perfect-props}(4)).
\end{proof}

\begin{cor} \label{cor:perfectns}
    The nullstellensatzians in $\CAlg_{\F_p}^{\sharp}$ are 
    discrete alg. closed fields.
\end{cor}



\begin{rmk} \label{rmk:eval-rke}
    \ref{thm:perfect-embed-discrete} and its proof
    together provide a technique for 
    evaluating the rke of an fpqc hypersheaf $F$
    on discrete perfect $\F_p$-algebras 
    to all perfect $\F_p$-algebras:
    
    Given a perfect $\F_p$-algebra $R$
    pick an hf hypercover $R_\bullet$ of $R$ by discrete perfect $\F_p$-algebras. 
    If $F(R_\bullet |_{\Delta_{\mathrm{inj}}})$ admits a limit, 
    then $(\iota^\heartsuit)_*F$ will exist pointwise at $R$ and agree with this limit.
\end{rmk}


\subsection{Unicity of perfect algebras} \label{sec:derived}

Up to this point we have confined ourselves to studying 
perfect commutative $\F_p$-algebras.
In this section we show that this in fact \emph{coincides} with 
perfect derived $\F_p$-algebras.
Derived algebras here being taken in the sense introduced by 
Brantner--Mathew in \cite{brantner2023deformation} and Raksit in \cite{raksit2020hochschild}.

\begin{cons}
    A derived ${\F_p}$-algebra is \mdef{perfect} 
    if its underlying commutative ${\F_p}$-algebra is perfect (in the sense of \ref{dfn:perfect}).
    This is captured by the pullback
    \[
    \begin{tikzcd}[sep=small]
	   {\rm{DAlg}_{\F_p}^\sharp} & {\rm{DAlg}_{\F_p}} \\
	   {\rm{CAlg}_{\F_p}^\sharp} & {\rm{CAlg}_{\F_p}}.
	   \arrow[from=1-1, to=1-2, hook, "\iota'"]
	   \arrow[from=1-1, to=2-1]
	   \arrow["\lrcorner"{anchor=center, pos=0.125}, draw=none, from=1-1, to=2-2]
	   \arrow[from=1-2, to=2-2, "un"]
	   \arrow[from=2-1, to=2-2, hook, "\iota"]
    \end{tikzcd}
    \] 
    As $un$ preserves colimits and limits \cite[Prop. 4.2.27]{raksit2020hochschild},
    the pullback above may be taken in either $\Pr^L$ or $\Pr^R$.
    Thereby, we learn that $\iota'$
    admits a left adjoint $(-)^{\sharp}$ (perfection) and 
    a right adjoint $(-)^{\flat}$ (tilting).
    We organize (perfect) derived $\F_p$-algebras into sites
    \[ 
    \DAff_{\F_p, \leq 0} \coloneqq \left( \DAlg_{\F_p, \leq 0}^{\op},\ \mathrm{hf} \right),\qquad
    \DAff^{\sharp}_{\F_p} \coloneqq \left( \DAlg_{\F_p}^{\sharp,\op},\ \mathrm{hf} \right).
    \]
\end{cons}

\begin{thm} \label{thm:unicity}
    The functor 
    $ \rm{DAlg}_{\F_p}^{\sharp} \rightarrow \rm{CAlg}_{\F_p}^{\sharp} $
    is an isomorphism.
\end{thm}

Our strategy for proving \ref{thm:unicity} is explicated by the next theorem.

\begin{thm} \label{thm:unicity-stacks}
    The maps of sites $\Aff_{\F_p}^\sharp \to \DAff^{\sharp}_{\F_p} \to \SAff_{\F_p}^{\sharp}$ induces isomorphisms 
    \[ \Shv(\Aff_{\F_p}^{\sharp}) \cong \Shv(\DAff_{\F_p}^{\sharp}) \cong \Shv(\SAff_{\F_p}^{\sharp}) \]
    on categories of accessible hypersheaves.
    Through \ref{cor:modules-sheaf} this creates an
    embedding of perfect commutative and perfect derived $\F_p$-algebras
    into $\Shv(\Aff_{\F_p}^{\sharp})^{\op}$.
    The essential images of these embeddings coincide proving \ref{thm:unicity}.
\end{thm}

As in \ref{thm:perfect-embed-discrete} the key input comes from control over the free perfect derived $\F_p$-algebra on a class in negative degree. Less perfectness this is provided by a calculation of Priddy:

\begin{rec} \cite[Thm.~4.1.1, Thm.~4.2.1]{priddy1973mod} \label{priddyfreealgebra}
    Let $\F_p [x_n]$ denote the free derived $\F_p$-algebra on a class in degree $n < 0$.
    The graded ring $\pi_*\F_p[x_n]$ is a free graded-commutative ring on classes in degrees $\leq n$.
\end{rec}

\begin{prop}[Toen \cite{toen2006affinestackschampsaffines}] \label{prop:Toen}
    The inclusion of sites $\Aff_{\F_p} \hookrightarrow \DAff_{\F_p, \leq 0}$ 
    induces an isomorphism on categories of accessible hypersheaves.
    Through \ref{cor:modules-sheaf} this creates an
    embedding of coconnective derived $\F_p$-algebras
    into $\Shv(\Aff_{\F_p})^{\op}$.
\end{prop}

\begin{proof}
    The proof is similar to \ref{thm:perfect-embed-discrete}.
    Combining  
    \ref{exm:hf-is-fpqc},
    \ref{exm:hf-cover-easy},
    \ref{lem:hf-props}(2) and
    \ref{priddyfreealgebra}
    we learn that the maps 
    $\{  \F_p[x_n] \to \F_p  \}_{n < 0}$ are 
    hf-covers verifying variant \ref{var:small-object} of \ref{lem:topos-basis}.
\end{proof}


\begin{lem} \label{lem:perfect-derived}
    The $\E_2$-algebra (in spaces) of endomorphisms of the identity on $\DAlg_{\F_p, \leq 0}$ 
    is isomorphic to $\N$ with generator the frobenius, $F$.
    
\end{lem}

\begin{proof}
    From \ref{prop:Toen} we obtain an isomorphism between the space of interest
    and the space of endomorphisms of the identity on affine $\F_p$-schemes.
    The identity functor on the latter category is represented by the affine line ring scheme $\A$.
    The space of ring scheme endomorphisms of $\A$ is isomorphic to the set of 
    polynomials $g(x) \in \F_p[x]$ that are 
    additive ($g(y+z) = g(y) + g(z)$), 
    multiplicative ($g(yz) = g(y) g(z)$) and
    unital ($g(1) = 1$).
    The only such polynomials are $x^{p^k}$---corresponding to the frobenius and its iterates.
\end{proof}

\begin{lem} \label{lem:colim-perfection}
    Given an $R \in \DAlg_{\F_p, \leq 0}$ we have $ \colim_F R \cong R^{\sharp} $.
\end{lem}

\begin{proof}
    For discrete $R$ this is clear.
    Now observe that the collection of $R$ for which this holds is closed under semi-cosimplicial limits.
    \ref{prop:Toen} witnesses every $R \in \DAlg_{\F_p, \leq 0}$ as such a limit of discrete rings.
\end{proof}

\begin{proof}[Proof (of \ref{thm:unicity-stacks}).]
    Let $\F_p [x_n]$ denote the free derived $\F_p$-algebra on a class in degree $n$.
    Consider the inclusion $\Aff_{\F_p}^\sharp \hookrightarrow \DAff^{\sharp}_{\F_p}$.
    Combining  
    \ref{exm:hf-is-fpqc},
    \ref{exm:hf-cover-easy},
    \ref{lem:hf-props}(2), 
    \ref{priddyfreealgebra} and
    \ref{lem:colim-perfection}
    we learn that the maps 
    $\{  \F_p[x_n]^{\sharp} \to \F_p  \}_{n < 0}$ are 
    hf-covers verifying variant \ref{var:small-object} of \ref{lem:topos-basis}.
    Together with \ref{thm:perfect-embed-discrete} this completes the proof of main claim.
    This yields a diagram 
    \[ \begin{tikzcd}[sep=small]
        \DAlg_{\F_p}^{\sharp} \ar[rr, hook] \ar[dr, "y_{0}"', hook] & &
        \CAlg_{\F_p}^{\sharp} \ar[dl, "y_{1}", hook] \\
        & \Shv(\Aff_{\F_p}^{\sharp})^{\mathrm{op}}.
    \end{tikzcd} \]
    
    What remains is showing that the horizontal arrow is essentially surjective.
    In fact, as this map is already known to be fully faithful and preserve colimits, we need only
    show that its image contains a family of generators.
    For this we calculate that
    \[ y_{0}( \F_p [x_n]^{\sharp} ) \cong \left( R \mapsto K(R, -n) \right) \cong y_{1}( \F_p \{x_n\}^{\sharp} ).\hfill\qedhere \]    
\end{proof}

\begin{cor} \label{cor:aut-pefect}
    The $\E_2$-algebra (in spaces) of endomorphisms of the identity on $\CAlg_{\F_p}^{\sharp}$ 
    is isomorphic to $\Z$ with generator the frobenius.
\end{cor}

\begin{proof}
    Compare \ref{lem:perfect-derived}.
    Using \ref{thm:unicity-stacks} we reduce to the discrete setting.
    This means analyzing ring scheme endomorphisms of the perfected affine line $\A^{\sharp}$.
    Concretely, we want additive, multiplicative, unital $g \in \F_p[x^{1/p^\infty}]$.
    The only such are $x^{p^k}$ for $k \in \Z$.
\end{proof}

While \ref{prop:Toen} provides each coconnective derived $\F_p$-algebra  
with an hf-hypercover by discrete $\F_p$-algebras,
it provides little control over the terms of this hypercover.
With some effort we can retain control over the first term
in the case where $\pi_0$ is a field.

\begin{lem} \label{lem:field-cover}
    Any $\pi_0$-iso $k \to A$ in $\DAlg_{k, \leq 0}$ admits a retraction.
\end{lem}

\begin{proof}
    Invoking \ref{prop:Toen}, pick a presentation of $A$ as a limit of a semicosimplicial diagram $R_\bullet$ of discrete $k$-algebras. 
    Let $B$ be the equalizer of $\delta_0, \delta_1 : R_0 \to R_1$.
    Observe that $A \to B$ is a $\pi_0$-iso with $B$ concentrated in $[-1,0]$.

    Proceeding by a small object argument it will now suffice to show that 
    for each $z \in \pi_{-1}(B)$ there exists a map $B \to B'$ sending $z$ to $0$.
    Let $C$ denote the free derived $B$-algebra with a nullhomotopy of $z$.
    Let $C'$ denote the free filtered derived $B$-algebra 
    (with filtration parameter $\tau$) equipped with a nullhomotopy of $\tau z$.
    On associated graded this is 
    the free graded derived $B$-algebra on a class $x$ in grading $1$, degree $0$. 
    We will analyze $\pi_*C$ using the $\tau$-localized $\tau$-bockstein sseq for $C'$.
    For degree reasons the powers of $x$ are the only classes in this sseq that support differentials.
    From the Kudo transgression theorem these differential are determined by
    $d_{p^k}(x^{p^k}) = (Q^0)^k(z)$.
    This allows us to conclude, inductively, that
    each page of this sseq has degree $0$ part $k[x^{p^j}]$ for some $j \geq 0$ and 
    degree $-1$ part a free module over this polynomial algebra.
    If the sseq does not degenerate at a finite page, 
    then we must have $\pi_0C \cong k$ so that we may take $B' \coloneqq C$.
    If the sseq degenerates at a finite page,
    then any choice of lift of surviving polynomial generator to $\pi_0C$
    becomes a witness that $\pi_0C \cong k[t]$ with $\pi_{-1}(C)$ a free $k[t]$-module. 
    Now take $B' \coloneqq C/t$.    
\end{proof}

\subsection{Spherical Witt vectors}
\label{sec:sphericalwitt}
\label{sec:Witt}

Up to this point we have confined ourselves to working relative over $\F_p$.
In this section we extend the unicity results of the previous section beyond this setting
by showing that bounded below perfect $\F_p$-algebras admit unique $p$-complete spherical lifts.
Our results extend those of \cite[Exm. 5.2.7]{lurie2018elliptic2} and \cite[Prop. 2.2]{chroNS}
on spherical Witt vectors.

\begin{rec} 
    There is a \mdef{spherical Witt vectors}---tilting adjunction
    \[ \W : \CAlg_{\F_p}^{\heartsuit,\sharp} \rightleftharpoons \CAlg_p^{\wedge} : (-)^{\flat_0} \]
    with 
    $\W$ fully faithful, 
    $\F_p \otimes \W \cong \mathrm{Id}$ and
    $\flat_0$ given by the classical tilt of $\pi_0(-)/p$.
\end{rec}



\begin{prop} \label{prop:w-plus-adj}
    In the bounded below setting, we refine the above to an adjunction
    \[ \W_+ : \CAlg_{\F_p,+}^{\sharp} \rightleftharpoons \CAlg_{p,+}^{\wedge} : (-)^{\flat_+} \]
    with 
    $\W_+$ fully faithful,
    $\F_p \otimes \W_+ \cong \mathrm{Id}$ and
    $(-)^{\flat_+} \cong (- \otimes \F_p)^{\flat}$.
\end{prop}

\begin{proof}
    Let $\CAlg_+^{\sharp} \hookrightarrow \CAlg_{p,+}^{\wedge}$ denote the full subcategory on those 
    $p$-complete bounded below commutative algebras $R$ with $\F_p \otimes R$ perfect.
    We prove the proposition via of
    the tilting adjunction of \ref{lem:perfect-props}
    and two subclaims:
    \begin{enumerate}
        \item For any $A \in \CAlg_+^\sharp$ and $B \in \CAlg_{p,+}^{\wedge}$ the map below is an isomorphism.
        \[ \Map_{\Ss}(A,B) \to \Map_{\Ss}(A, B \otimes \F_p) \]
        \item The map $- \otimes \F_p : \CAlg_+^\sharp \to \CAlg_{\F_p, +}^{\sharp}$ admits an inverse $\W_+$.
    \end{enumerate}
    
    Recall that the adjunction 
    $ \CAlg_{p,+}^{\wedge} \rightleftharpoons \CAlg_{\F_p,+}$ is comonadic (cf. \cite{bousfield1979localization}).
    
    \textit{Part (1).}
    Comonadicity allows us to reduce to the case where $B$ carries an $\F_p$-algebra structure.
    Using this structure construct the diagram below where 
    $C \coloneqq A \otimes \F_p$, 
    $D \coloneqq \F_p \otimes \F_p$,
    the vertical composites are isomorphisms
    and the bottom right vertical map is an isomorphism by \ref{lem:tilt-coproduct} (below). 
    \[ \begin{tikzcd}[row sep=small]
        \Map(A, B) \ar[r, "\cong"] \ar[d] &
        \Map_{\F_p}(C, B) \ar[d] &
        \Map_{\F_p}(C, B^{\flat}) \ar[d] \ar[l, "\cong"'] \\
        \Map(A, B \otimes \F_p) \ar[r, "\cong"] \ar[d] &
        \Map_{\F_p}(C, B \otimes_{\F_p} D )  \ar[d] &
        \Map_{\F_p}(C, (B \otimes_{\F_p} D)^{\flat}) \ar[l, "\cong"'] \ar[d, "\cong"] \\
        \Map(A, B) \ar[r, "\cong"]  &
        \Map_{\F_p}(C, B \otimes_{\F_p} \F_p)  &
        \Map_{\F_p}(C, (B \otimes_{\F_p} \F_p)^{\flat}) \ar[l, "\cong"'] 
    \end{tikzcd} \]
    
    \textit{Part (2).}
    In view of (1) it will suffice to argue that 
    $ - \otimes \F_p : \CAlg_+^{\sharp} \to \CAlg_{\F_p,+}^{\sharp} $
    admits a section. 
    On discrete perfect $\F_p$-algebras $\W$ provides such a section. 
    We proceed by rke.  
    Using the comonadicity from above we observe that
    $\W_+ \coloneqq (\iota^{\heartsuit})_* \W$ will exist and satisfy $\F_p \otimes \W_+ \cong \mathrm{Id}$
    as soon as
    \begin{enumerate}
        \item[(a)] $\mathrm{Id} \circ \iota^{\heartsuit} \cong \iota^{\heartsuit}$ witnesses $\mathrm{Id}$ as a right-kan-extension of $\iota^{\heartsuit}$ along $\iota^{\heartsuit}$ and
        \item[(b)] $- \otimes \F_p^{\otimes \bullet}$ preserves the limits exhibiting (a).
    \end{enumerate}
    Both (a) and (b) may be extracted from \ref{lem:rke-identify} (below).
\end{proof}

\begin{lem} \label{lem:tilt-coproduct}
    Let $A,\ f : B_0 \to B_1 \in \CAlg_{\F_p}$.
    If $A$ is bounded below and 
    $\pi_n(f)$ is an isomorphism for $n \leq 0$,
    then $(A \otimes f)^\flat$ is an isomorphism.
\end{lem}

\begin{proof}
    Extend $\pi_\star(f)$ to an exact sequence
    $ 0 \to K \to \pi_\star B_0 \xrightarrow{\pi_\star(f)} \pi_\star B_1 \to C \to 0 $
    in $\DLM$.
    In view of \ref{homotopygroupsoftilt} it will suffice to argue that for $T \in \{C,K\}$
    \[ \cdots \xrightarrow{Q^0} (\pi_\star A \otimes T)_k \xrightarrow{Q^0} (\pi_\star A \otimes T)_k \]
    is isomorphic to the zero pro-object for each $k \in \Z$.
    Expanding $(\pi_\star A \otimes T)_k$ as 
    \[ (\pi_\ell A \otimes T_{k-\ell}) \oplus \cdots \oplus (\pi_{k-1} A \otimes T_1) \]
    we observe that instability and the Cartan formula imply $Q^0$ acts via 
    a strictly upper triangular matrix (i.e. nilpotently).    
\end{proof}

\begin{lem} \label{lem:rke-identify}
    Any $ \CAlg_{\F_p,+}^{\sharp} \xrightarrow{u} \Mod_{\F_p,+}^{\leq 0} \xrightarrow{g} \mathcal{E}$
    in $\mathrm{Cat}^{\mathrm{Lex}}$
    is the rke'd from $\iota^{\heartsuit}$.
\end{lem}

\begin{proof}
    Let $R$ be an $\ell$-connective perfect $\F_p$-algebra.
    Let $R_\bullet$ be an hf hypercover of $R$ by discrete perfect $\F_p$-algebras 
    (note that these exist by the proof of \ref{thm:perfect-embed-discrete}).
    As $u$ is an hf hypersheaf by \ref{cor:modules-sheaf}, $u R_\bullet$ is a limit diagram.
    Let $uR\langle k \rangle$ denote the limit of $uR_\bullet$ over $\Delta_{\mathrm{inj}}^{\leq k}$. 
    As each term of $uR_\bullet$ lies in the heart we learn that for $k+\ell > 0$ the fiber
    $ F_k \to uR \to uR\langle k \rangle $ is concentrated entirely in degree $-(k+1)$. 
    In particular, the tower $ uR \langle \bullet \rangle $ is quickly converging in the sense of 
    \cite[Dfn. 2.16]{clausen2017descentalgebraicktheoryconjecture}.
    In this way we learn that $guR_\bullet$ is a limit diagram.
    In view of \ref{rmk:eval-rke} the proof is now complete.
\end{proof}

\begin{rmk} \label{rmk:flat-plus}
    The functor $(-)^{\flat_+}$ witnesses that the tilt of a bounded below, commutative $\F_p$-algebra $R$ independent of the choice of $\F_p$-algebra structure.
\end{rmk}

\begin{rmk}
    Beyond the bounded below case it is unclear in what generality spherical lifts exist.
    As a positive example we have $\bb{S} \otimes _{C^{*}(S^3;\bb{S})} \bb{S}$ as a spherical lift of $C^{*}(\Omega S^3; \bb{F}_p)$.
\end{rmk}



    



\section{Spaces via cochains}
\label{sec:cochains}

Fix once and for the remainder of the section a perfect field $k$ of characteristic $p$. As observed in \ref{ex:cochainperf}, the $k$-cochains on a space $X$ are a perfect $k$-algebra. 

In this section, leveraging our understanding of perfect $k$-algebras, we give simplified proofs of two deep theorems about $p$-complete spaces: Mandell's theorem \cite{mandell2001algebras} and Miller's theorem \cite{miller1984sullivan} (nee the Sullivan conjecture). 
Our proof of Mandell's theorem is close in spirit to the original, but avoids providing an explicit presentation of $\cochain{K(\Z/p,n)}{k}$ 
in favor of an induction on $n$ whose base case ($n=0$)
is essentially trivial.
Building on our work in \cref{sec:sphericalwitt} we then prove a generic version of Yuan's variant \cite[Cor. 7.6.1]{yuan2023integral} of Mandell's theorem allowing $k$-cochains to be replaced by $R$-cochains for appropriate $R$.
Our proof of Miller's theorem is really a proof that
the Segal conjecture implies the Sullivan conjecture
(cf. \cite{carlsson1991equivariant}).
Our argument is a remarkably simple 
derivation given an observation developed in 
\cref{sec:fp-tilt}:
tilting provides an inverse Tate construction.

\subsection{Mandell's theorem} 
\label{sec:MandellThm}
\label{secpadichomotopy}


    

\begin{thm}[Mandell] \label{thm:mandell}
    For $k$ alg. closed, 
    the $k$-cochains functor 
    \[ \cochain{-}{k}: \Spc \to \rm{CAlg}^{op}_k \]
    is fully faithful whenever the target is 
    nilpotent and $p$-completely finite type.
\end{thm}

\begin{proof}
    Reduce to the case where the source is a point.
    Let $\A$ denote the collection of 
    nilpotent and $p$-completely finite type spaces
    for which the unit map 
    \[ \eta[X] : X \to \Map_k(C^*(X;k), k) \]
    is an isomorphism.
    We will prove the theorem by showing that $\A$ contains 
    all $p$-complete nilpotent spaces of finite type.

    \begin{itemize}
        \item[(a)] $\A$ is closed under finite coproducts.
        \item[(b)] $\A$ is closed under those limits preserved by $C^*(-;k)$.
        \item[(c)] $C^*(-;k)$ preserves fiber sequences
            $ F \to E \to B $ with $B$ connected, finite type and $\pi_1(B,b)$ a finite $p$-group.
            Via (b) this implies that if $*,B, E \in \A$, then $F \in \A$. 
        \item[(d)] $\A$ is closed under delooping.
        \item[(e)] $\A$ contains all $p$-complete nilpotent spaces of finite type.
    \end{itemize}
    
    \textit{Part (a).}
    As $C^*(-;k)$ preserves colimits this follows from 
    $\CAlg_k^{\op}$ being disjunctive.
    
    \textit{Part (b).}
    Suppose we are given a limit diagram of spaces that $C^*(-;k)$ takes to a colimit diagram of perfect commutative $k$-algebras.
    As $\Map_k(-,k)$ takes colimits to limits 
    our hypotheses imply that the comparison map of interest is a limit of isomorphisms and therefore itself an isomorphism.


    \textit{Part (c).}    
    Given $X \to B$,
    let $R_X \coloneqq \cochain{E}{k} \otimes_{\cochain{B}{k}} \cochain{X}{k}$ and 
    let $ c_X : R_X \to \lim_{x \in X} R_x  $ denote the coassembly map. 
    Given $b \in B$, 
    let $F_b \coloneqq E \times_B \{b\}$,
    let $G_b \coloneqq \pi_1(B,b)$ and
    let $ \theta_b : R_b \to \cochain{F_b}{k} $
    denote the assembly map of interest. 
    The maps $\theta_b$ and $c_B$ are related by the diagram below.
    \begin{equation} \label{eqn:coass}
    \begin{tikzcd}
    {\cochain{E}{k}} \ar[d, equal] \ar[r,"\cong"] &
    {R_B} \ar[r, "c_B"] \ar[d, "\cong"] &
    {\lim_{b \in B} R_b } \ar[d, "\lim_{b \in B} \theta_b"] \\
    {\cochain{E}{k}} \ar[r, "\cong"] &
    {\cochain{\colim_{b\in B} F_b}{k}} \ar[r, "\cong"] &   
    {\lim_{b \in B} \cochain{F_b}{k}}   
    \end{tikzcd}
    \end{equation} 

    Using (\ref{eqn:coass}) we will reduce (c) to verifying that $c_B$ is an isomorphism.
    Specifically, we will show that 
    $\fib(\theta_b) = 0$ if and only if $\lim_{b \in B} \fib(\theta_b) = 0$.
    From \ref{exm:hf-cover-easy}, $\cochain{B}{k} \to \cochain{E}{k}$ is an hf cover
    and $\cochain{E}{k}$ therefore has tor-amplitude $\leq 0$ 
    as a right $\cochain{B}{k}$-module.
    This implies $\fib(\theta_b)$ is bounded above.
    If $\fib(\theta_b)$ is $(n+1)$-coconnective, then
    $ \pi_n\left( \lim_{b \in B} \fib(\theta_b) \right) \cong \pi_n(\fib(\theta_b))^{G_b} $.
    Varying $n$ and recalling that $G_b$ is a finite $p$-group 
    (reps of $p$-groups in char. $p$ are unipotent)
    completes the reduction.
        
    Let $\mathbb{B}_n$ denote the collection of $X \to B$ for which $c_X$ is $n$-coconnective. 
    It is clear that $\mathbb{B}_n$ is closed under finite colimits, retracts and contains all $* \to B$.
    What is less clear (but also true) is that
    given $f : X \to Y \in \Spc_{/B}$ 
    with $H_i(f;k)$ an isomorphism for $i \leq -n$,
    we have $X \in \mathbb{B}_n$ if and only if $Y \in \mathbb{B}_n$.
    For this we observe that on the one hand,
    our tor-amplitude bound on $\cochain{E}{k}$ implies 
    $s(c_f)$ is an isomorphism in degrees $\geq n$.
    While on the other hand,
    $t(c_f)$ may be analyzed by 
    filtering $R_x$ by a refinement of the postnikov tower 
    chosen so that the associated graded has trivial action 
    (recall that each $G_x$ is a finite $p$-group),
    reducing to this associated graded,
    then evaluating directly.
    The upshot of this again being that the map is an isomorphism in degrees $\geq n$.
    Now invoke $k$-homology skeletal approximations
    to witness that $B \in \mathbb{B}_n$ for all $n$.
    
    \textit{Part (d).}
    Let $(X,x)$ be a pointed connected space with $\Omega_x X \in \A$.
    Consider a pushout 
    \[ \begin{tikzcd}
        C^*(X;k) \ar[r, "\nu"] \ar[d, "x"] & 
        k \ar[d, "x'"] \\
        k \ar[r, "\nu'"] &
        P_{x,\nu} \pushout
    \end{tikzcd} \]
    where $\nu$ is arbitrary.
    From our hypothesis that $\Omega_x X \in \A$ and (c) we obtain isomorphisms 
    $P_{x,x} \cong C^*(\Omega_x X; k)$ and $\Omega_x \eta[X] \cong \eta[\Omega_x X]$.
    Thus (d) is reduced to showing that 
    $\Map_{k}(P_{x,\nu}, k)$ is non-empty for all $\nu$.
    
    In this, we start by using \ref{exm:hf-cover-easy}, \ref{lem:hf-props} 
    to conclude that $x$ and $x'$ are hf-covers.
    In particular, $P_{x,\nu} \neq 0$.
    Next, we calculate that
    \begin{align*}
        P_{x,\nu}\ {}^{x'}\!\!\otimes_k^{x'} P_{x,\nu}
        &\cong (k\ {}^{x}\!\!\otimes_{C^*(X;k)}^{x} k) \otimes_{C^*(X;k)}^{\nu} k \\
        &\cong P_{x,x} \otimes_{C^*(X;k)}^{\nu} k \\
        &\cong C^*(\Omega_x X;k) \otimes_{k}^{\nu'} P_{x,\nu}.
    \end{align*}
    from which we learn that $\pi_0P_{x,\nu}$ is finite over $k$
    of degree $|\pi_0 \Omega_x X|$
    and therefore isomorphic to $k^{\pi_0 \Omega_x X}$
    (here we use that $\pi_0P_{x,\nu}$ is perfect and $k$ is alg. closed).    
    Now pick a $z \in \pi_0P_{x,\nu}$
    such that $\pi_0 P_{x,\nu}[z^{-1}] \cong k$
    so that by \ref{lem:field-cover}
    there exists a map $P_{x,\nu}[z^{-1}] \to k$.
    
    
    \textit{Part (e).}
    Applying (b) to the empty limit we conclude that $\A$ contains the point.
    Applying (a) to the point we conclude $\A$ contains all finite sets.
    Applying (d) we conclude $\A$ contains $B\Z/p^j$.
    Observe that the limit diagram 
    \[ B\Z_p \to \cdots \to B\Z/p^2 \to B\Z/p \]
    is preserved by $C^*(-;k)$.
    Applying (b) we learn that $\A$ contains $B\Z_p$. 
    
    From the theory of Moore--Postnikov towers 
    we know that the collection of spaces 
    generated from $BC_p$ and $B\Z_p$ 
    under the operations from (c) and (d)
    is the collection of all
    connected, bounded, nilpotent, $p$-completely finite type spaces.

    The boundedness hypotheses can be removed by 
    observing that limits of Postnikov towers commute with taking cohomology and applying (b). 
    The connectedness hypothesis can be relaxed to finiteness of $\pi_0$ using (a). 
\end{proof}

\begin{cor}[Generic Mandell--Yuan] 
\label{cor:yuan-model} 
    Let $R \in \CAlg_{p,+}^{\wedge}$ with $k \cong R^{\flat_+} $ alg. closed.
    The $R$-cochains functor is fully faithful
    at nilpotent, $p$-complete, $\F_p$-small\footnote{%
    Meaning: $X$ lies in the full subcategory of spaces generated from $*$ under finite colimits, retracts and $\F_p$-homology isomorphisms.}
    target spaces.
\end{cor}

\begin{proof}
    Reduce to the case where the source is a point.
    Using our hypotheses on $X$ and $R$, 
    the adjunction from \ref{prop:w-plus-adj} and
    the induced unit map $\W(k) \to R$
    we calculate that 
    \begin{align*}
    \Map_{{R}}(\cochain{X}{R}, R) 
    &\cong \Map_{\W(k)}(\cochain{X}{\W(k)}, R) \\
    &\cong \Map_{\W(k)}(\W_+\cochain{X}{k}, R) \\
    &\cong \Map_{k}(\cochain{X}{k}, k)    
    \end{align*}
    reducing the corollary to \ref{thm:mandell}.
\end{proof}

\begin{exm}
    Let $k = \overline{k}$.
    The $p$-completions of each of the commutative algebras
    $\W(k)$, 
    $\mathrm{MU} \otimes \W(k)$,
    $\mathrm{ku} \otimes \W(k)$ and
    $\Z \otimes \W(k)$
    satisfy the hypotheses of \ref{cor:yuan-model}.
\end{exm}

\subsection{Fixed points via tilting}
\label{sec:fp-tilt}




In this section we briefly revisit a classic question:
Given a $C_p$-space $X$ what 
information about $X^{C_p}$
can we recover from $X^e$ with its $C_p$-action.
Our starting point will be 
the fundamental result of this type: Borel's localization theorem.

\begin{rec} \label{rec:Borel}
    Recall that 
    $\cochain{BC_p}{k} \cong k^{BC_p} \cong k^{hC_p}$ 
    has homotopy groups 
    $R \coloneqq k[e,t]/(e^2 +(1-\ptwo{1})t)$
    where $|e|=-1$ and $|t|=-2$ while 
    $\cochain{BC_p}{k}[t^{-1}] \cong k^{tC_p}$
    is the localization of this at $t$.
    
    Let $X$ be an $\F_p$-small\footnote{%
    Meaning: $X$ lies in the full subcategory of $C_p$-spaces generated from $*$ and $C_p$ under finite colimits, retracts and $\F_p$-homology isomorphisms.} 
    $C_p$-space.
    The fixed point inclusion $X^{C_p} \to X^e$ induces isos
    \begin{align} \label{eqn:localization}
        \cochain{X^{C_p}}{k}^{tC_p} &\cong \cochain{X^e}{k}^{tC_p}, \\
       H^*(X^{C_p};k) \otimes_k R[t^{-1}] &\cong H^*(X_{hC_p};k)[t^{-1}].
    \end{align}
\end{rec} 

\begin{proof}
    Both terms of interest are $\F_p$-homology isomorphism invariant and preserve finite colimits, retracts.
    We may therefore reduce to 
    the cases $X=*$ and $X=C_p$.
    For $X=*$ the fixed point inclusion itself is an isomorphism.
    For $X=C_p$ the fixed points are empty
    and $\cochain{C_p}{k}^{tC_p} = 0$.
\end{proof}

Examining (\ref{eqn:localization}), the salient question is 
how much information about $X^{C_p}$ does this iso provide?
The main objective of this section is 
to demonstrate that tilting delocalizes (\ref{eqn:localization})
into a formula for $\cochain{X^{C_p}}{k}$.
In view of Mandell's theorem this leads to the striking conclusion that $\cochain{X^e}{k}^{tC_p}$
records essentially complete information on $X^{C_p}$.



\begin{prop} \label{prop:tilt-undo-k}
    Let $A \in \CAlg_{\F_p}^{\sharp, \heartsuit}$.
    The natural map 
    $A^{hC_p} \to (A^{tC_p})^{\flat}$
    is iso.
\end{prop}

\begin{proof}
    The map $A^{hC_p} \to A^{tC_p}$ induces an iso on $\pi_n$ for $n \leq 0$.
    In particular, $Q^0$ already acts by an iso
    on the non-positive homotopy groups of $k^{tC_p}$.
    Now apply \ref{homotopygroupsoftilt}, \ref{rmk:tilt-simplify}.
\end{proof}

\begin{cnstr}
    Let $\Phi$ denote the functor
    \[ \Phi : ( \CAlg_{\F_p}^{\sharp} )^{BC_p} \to \CAlg_{\F_p}^{\sharp}, \qquad \Phi A \coloneqq \F_p \otimes_{\cochain{BC_p}{\F_p}} (A^{tC_p})^{\flat} \]
    where the $\cochain{BC_p}{\F_p}$-algebra structure on the tilt comes from \ref{prop:tilt-undo-k}.
\end{cnstr}


\begin{thm} \label{tiltgenuinefixedpoints}
    Let $X$ be an $\F_p$-small $C_p$-space.
    We construct a natural isomorphism of perfect $k$-algebras
    \begin{equation} \label{eqn:tilt-localization}
        \cochain{X^{C_p}}{k} \cong \Phi C^*(X^e;k).     
    \end{equation}
\end{thm}

\begin{proof}
    In view of (\ref{eqn:localization})
    it will suffice to show that the coassembly map
    $\cochain{Y}{k} \to \Phi \cochain{Y}{k}$
    associated to $\Phi k \cong k$ (see \ref{prop:tilt-undo-k}) 
    is an iso for $Y$ an $\F_p$-small space.
    Both source and target of this map are $\F_p$-homology isomorphism invariant and preserve finite colimits, retracts.
    We thereby reduce to the case $Y=*$
    where we recover the iso 
    used in the construction of the coassembly map.
\end{proof}



\begin{rmk}
    Combining \ref{tiltgenuinefixedpoints} and \ref{thm:mandell}, we learn that the 
    $p$-completion of the $C_p$-fixed points of $X$ can be recovered from the
    underlying Borel $C_p$-space of $X$. This can already be viewed as a weak form of the Sullivan conjecture.
\end{rmk}



\subsection{Segal implies Sullivan}


Using the whole range of things developed in this paper we now give a short proof that the Segal conjecture for $C_p$ (as proved in \cite{lin1980conjectures, lin1980calculation, gunawardena1980segal}) implies the Sullivan conjecture (as proved in \cite{miller1984sullivan}, \cite{lannes1992espaces}, \cite{carlsson1991equivariant}).

\begin{thm} \label{thm:sullivan}
    Let $X$ be a nilpotent, $p$-complete, $\F_p$-small $C_p$-space.
    Then $X$ is Borel-complete in the sense that the natural map 
    $ X^{C_p} \to (X^e)^{hC_p} $
    is an isomorphism.
\end{thm}

\begin{proof}
    Let $k = \overline{k}$.
    We prove the theorem via two subclaims:
    \begin{enumerate}
        \item Let $q : BC_p \to *$. The diagonal map 
        $\Map_{\Spc}(q, X^{C_p})$ is iso.
        \item Let $\iota : X^{C_p} \to X^e$ denote the fixed points inclusion. $\iota^{hC_p}$ is iso.
    \end{enumerate}
    \medskip
    
    Part (1).
    Using the Segal conjecture\footnote{%
    More precisely, the extension of the Segal conjecture proven in %
    \cite[Prop. 6.6]{yuan2023integral}, \cite[Thm. 3.3]{burklund2024notesegalconjecturelarge}.}
    and \ref{homotopygroupsoftilt} we
    observe that $\cochain{BC_p}{\W(k)}$ is connective
    and calculate that
    \begin{equation} \label{eqn:tilt-switt-BCp}
        \cochain{BC_p}{\W(k)}^{\flat_+} \cong \left( \W(k)^{BC_p} \right)^{\flat_+} \cong (k \otimes A(C_p))^{\flat} \cong (k[x]/x^2)^{\flat} \cong k
    \end{equation}
    where $A(C_p)$ is the Burnside ring of $C_p$.
    Using \ref{cor:yuan-model} 
    and \ref{prop:w-plus-adj}
    we find that
    \begin{align*}
        \Map_{\Spc}(q,X^{C_p})
        &\cong \Map_{\W(k)}(\cochain{X^{C_p}}{\W(k)}, \cochain{q}{\W(k)}) \\
        &\cong \Map_{\W(k)}(\W_+\cochain{X^{C_p}}{k}, \cochain{q}{\W(k)}) \\
        &\cong \Map_{k}(\cochain{X^{C_p}}{k}, \cochain{q}{\W(k)}^{\flat_+}) 
    \end{align*}
    with the final term being an iso by     (\ref{eqn:tilt-switt-BCp}).
    
    Part (2).
    Using 
    straightening--unstraightening,
    \ref{thm:mandell} and
    \ref{prop:tilt-undo-k}
    we find that
    \begin{align*}
        \iota^{hC_p}
        &\cong \Map_{\Spc_{/BC_p}}(BC_p, \iota_{hC_p}) \\
        &\cong \Map_{\cochain{BC_p}{k}}(\cochain{\iota_{hC_p}}{k},\cochain{BC_p}{k}) \\
        &\cong \Map_{k^{hC_p}}(\cochain{\iota}{k}^{hC_p}, (k^{tC_p})^{\flat}) \\
        &\cong \Map_{k^{hC_p}}(\cochain{\iota}{k}^{hC_p}, k^{tC_p}) \\
        &\cong \Map_{k^{hC_p}}(\cochain{\iota}{k}^{tC_p}, k^{tC_p})
    \end{align*}
    with the final term being an iso by (\ref{eqn:localization}).
\end{proof}

\printbibliography

\end{document}